\documentclass[11pt,a4paper]{article}

\usepackage[utf8]{inputenc}
\usepackage{amsmath, amsthm, amssymb}
\usepackage{graphicx}
\numberwithin{equation}{section}

\begin{document}

\newtheorem{theorem}{Theorem}
\newtheorem{lemma}{Lemma}
\newtheorem{definition}{Definition}
\newcommand{\ob}{\textbf}

\title{A regularization for the transport equations using spatial-averaging}  

 \author{John Villavert\thanks{Graduate Student, Department of
         Applied Mathematics, University of Colorado,
         Boulder, Colorado, 80309, US.}\and Kamran Mohseni\thanks{Associate Professor of
        Aerospace Engineering Sciences; Affiliated faculty
        in the Applied Mathematics Department, University of
        Colorado, Boulder, Colorado, 80309, US.}
        }
\date{}
\maketitle

\begin{abstract}
This paper examines an averaging technique applied to the transport equations as an alternative to vanishing viscosity. Such techniques have been shown to be valid shock-regularizations of the Burgers equation and the Euler equations, but has yet to be applied to the similar transport equations. However, for this system, the classical notion of weak solutions is not always sufficient thus a more general notion of a distribution solution containing Dirac-delta functions must be introduced. Moreover, the distribution solution to the Riemann problem is known to be the weak-$*$ limit of the viscous perturbed transport equations as viscosity vanishes. In comparison to the classical method of vanishing viscosity, the Riemann problem is examined for the averaged transport equation and it is shown that the same delta-shock distribution solution is captured as filtering vanishes. Both mathematical and physical motivation are provided for the averaging techniques considered including an existence and uniqueness result of solutions within the class of smooth initial conditions.
\end{abstract}


\section{Introduction} 
Consider the transport (or transportation) equations in multiple dimensions in which the continuity equation is adjoined with the inviscid Burgers equation
\begin{equation}
  \left\{\begin{array}{r@{\;=\;}l@{\quad}l}
      	\rho_{t} + \nabla \cdot (\rho \ob{u}) & 0 \\
 	 			\ob{u}_{t} + \ob{u} \cdot \nabla \ob{u} & 0 \\
 	 			(\rho(\ob{x},0),\ob{u}(x,0)) & (\rho_{0}(x),\ob{u}_{0}(x)).
        \end{array}
 \right.
\end{equation}
The transport equations provide a simple system of two conservation laws for which no classical weak solutions may exist for some initial data, so a distribution solution consisting of Dirac delta functions or a delta-shock solution must be introduced \cite{Sheng:99a, Tan:94a}. On the other hand, the transport equations model the dynamics of particles that adhere to one another upon collision and has been studied as a simple cosmological model for describing the nonlinear formation of large-scale structures in the Universe \cite{Frisch:01a}. For example, $\ob{u}$ may represent the flow field carrying dust particles with density $\rho$, and the delta-shock wave represents a concentration of dust on a shock which attract the dust \cite{Shen:09a}. For the Riemann problem, the classical regularization of adding dissipation such as the adhesion model \cite{Frisch:01a,Gurbatov:91a},
\begin{equation}
  \left\{\begin{array}{r@{\;=\;}l@{\quad}l}
      	\rho_{t} + \nabla \cdot (\rho \ob{u}) & 0 \\
 	 			\ob{u}_{t} + \ob{u} \cdot \nabla \ob{u} & \nu \Delta \ob{u},
        \end{array}
 \right.
\end{equation}
have been shown to recover the delta-shock solution as $\nu \rightarrow 0$. Similarly, K.T. Joseph \cite{Joseph:93a} considered the Riemann problem for the following regularization
\begin{equation}
			\begin{array}{r@{\;=\;}l@{\quad}l}
      	\rho_{t} + (\rho u)_{x} & \epsilon u_{xx} \\
 	 			u_{t} + (\frac{u^{2}}{2})_{x} & \epsilon u_{xx}
      \end{array}
\end{equation}
and proved that the limit of the corresponding solutions $(\rho^{\epsilon},u^{\epsilon})$ may contain delta-measures as $\epsilon \rightarrow 0$. The authors in \cite{Tan:94a} examined a similar system in one dimension
\begin{equation}
  \begin{array}{r@{\;=\;}l@{\quad}l}
      	\rho_{t} + (\rho u)_{x} & 0 \\
 	 			u_{t} + (u^{2})_{x} & 0
  \end{array}
\end{equation}
in order to mathematically justify delta-shock waves. They, too, introduce viscous perturbations. One case involves adding the typical Laplacian term
\begin{equation}
  \begin{array}{r@{\;=\;}l@{\quad}l}
      	\rho_{t} + (\rho u)_{x} & 0 \\
 	 			u_{t} + (u^{2})_{x} & \epsilon u_{xx}
  \end{array}
\end{equation}
while the other introduces time-dependent viscosity
\begin{equation}
  \begin{array}{r@{\;=\;}l@{\quad}l}
      	\rho_{t} + (\rho u)_{x} & 0 \\
 	 			u_{t} + (u^{2})_{x} & \epsilon tu_{xx}.
        \end{array}
\end{equation}
They construct solutions to these perturbed systems and show that the delta-shock waves are the weak-$*$ limits of the the solutions to the viscous systems as viscosity vanishes.

The principle focus of this paper is to introduce an alternative regularization of the transport equations in which the delta-shock solution is captured similar to what vanishing viscosity achieves, yet without largely altering the initial structure of the equations. This particular regularization will utilize a spatial averaging of the nonlinear terms by convoluting them by some averaging kernel or low-pass filter rather than adding viscous terms. The notion of averaging the transport equations directly stems from the successful regularization of the Burgers equation via averaging of the convective velocity \cite{Mohseni:06l},  
\begin{equation}
  \left\{\begin{array}{r@{\;=\;}l@{\quad}l}
 	 			u_{t} + \overline{u} u_{x} & 0 \\
 	 			\overline{u} & g^{\alpha}*u \\
 	 			g^{\alpha} & \frac{1}{\alpha} g(\frac{x}{\alpha}).
        \end{array}
 \right.
\end{equation}
Their numerical study suggests that this filtering of the nonlinear term attenuates shock formation without adding viscosity or dispersion. Mohseni and Norgard \cite{Mohseni:08w,Mohseni:09h} as well as  Bhat and Fetecau \cite{BhatHS:06a} later studied this regularization in one and multiple dimensions. It is shown that spatial averaging of the convective velocity properly regularizes the Burgers equation while numerical simulations suggest the entropy solution is captured as the filtering parameter, $\alpha$, vanishes. In fact, the entropy solution is captured for a class of `bell-shaped' $C^{1}$ initial conditions \cite{Mohseni:09h}. 

It should be noted that this was not the first implementation of filtered convective velocities. In 1934, Jean Leray proposed using a filtered convective velocity in the Navier-Stokes equations \cite{LerayJ:34a}. This study has influenced the investigation of the Leray-$\alpha$ models of turbulence \cite{Holm:04a} as well us other models of turbulent flows. The Lagrangian Averaged Navier-Stokes-$\alpha$ (LANS-$\alpha$) uses similar filtering to successfully model some turbulent incompressible flows \cite{Marsden:98b,Foias:01a,Mohseni:03a}. Although we do not use Lagrangian type averaging in this manuscript, it has ultimately inspired the concept of spatial filters on equations of compressible flows as a potential model for shock formation and turbulence.

The spatial averaging in the CFB equation has also been extended to both the compressible homentropic Euler and full Euler equations in \cite{BhatHS:09b,Mohseni:09w,Mohseni:09j}. For instance, Bhat and Fetecau \cite{BhatHS:09b} applied a Leray-type averaging with the Helmholtz filter $g$ to obtain the equations
\begin{equation}
  \begin{array}{r@{\;=\;}l@{\quad}l}
      	\rho_{t} + \overline{u}\rho_{x} + \rho u_{x} & 0 \\
 	u_{t} + \overline{u} u_{x} + \frac{\kappa \rho_{\gamma}}{\rho} & 0,
  \end{array}
\end{equation}
where
\[ g(x) = \frac{1}{2\alpha} exp\left(-\frac{|x|}{\alpha}\right). \]
For the isothermal case, $\gamma = 1$, the authors show the averaging prevents the crossing of characteristics and prove the unique existence of global smooth solutions. They further show that the solutions will converge strongly to a weak solution of the homentropic Euler equations as $\alpha \rightarrow 0$, however, they conclude from numerical simulations that this is not the physically relevant solution (Note: for $\gamma \neq 1$, this Leray-type averaging does not prevent characteristics from crossing). 

A promising averaging technique introduced in \cite{Mohseni:09w}, called the observable divergence method, was applied to both the homentropic Euler and full Euler equations in one dimension \cite{Mohseni:09j}. Mohseni \cite{Mohseni:09w} derived these equations from basic principles. Unlike the aforementioned Leray-type averaging, the observable divergence method appeared to numerically capture the entropy solutions for the shock-tube problem and the Shu-Oscher problem. This particular type of averaging will receive full attention within this manuscript since it appears to be a valid shock-regularization for the equations of gas dynamics. The importance of this regularization is its potential to model small scale behavior capable of capturing shocks and (possibly) turbulence in a single technique, thus reducing the computational cost demanded by many problems in engineering flows. 

Influenced by such work, we shall adopt spatial averaging for the transport equations by applying the method of observable divergence to the continuity equation and Burgers equation then obtain analogous results to the method of vanishing viscosity. This manuscript will be organized as follows. Section \ref{Preliminaries} will review basic material on the transport equations, delta-shock waves, and the Riemann problem. Sections \ref{Filters} and \ref{Averaging} will discuss the class of averaging kernels in consideration and will introduce and briefly motivate the method of observable divergence. A derivation of the observable transport equation is provided. Section \ref{Observable equations} examines the observable transport equations with an existence and uniqueness result in the class of smooth initial conditions. The one-dimensional case is considered with a discontinuous initial condition consisting of a single decreasing jump since this provides a Dirac delta distribution solution for the original transport equations. It is shown that the regularized equations share the same Dirac delta distribution solution.

\section{Preliminaries}\label{Preliminaries}
Let us first review several fundamental ideas and results for the transport equations. This section reviews the one-dimensional transport equations and its distribution solution for Riemann initial data.

\subsection{The one-dimensional initial value problem}
For now, let us restrict our attention to the one-dimensional initial value problem
\begin{align}
	\rho_{t} + (\rho u)_{x} = 0 \label{onedim1} \\ 
	u_{t} + (\frac{u^{2}}{2})_{x} = 0 \label{onedim2} \\ 
	(\rho(x,0),u(x,0)) = (\rho_{0}(x),u_{0}(x)) \label{onedim3}	
\end{align}

This gives us a simple example of a system of a nonlinear nonstrictly hyperbolic PDE in which classical weak solutions may not exist. Instead, one must formulate a broader notion of a distribution solution and include a new type of nonlinear hyperbolic wave, a delta-shock wave. To illustrate this idea, assume that the initial conditions are smooth and $u_{0}'(s) < 0$ for some point s. Then the characteristic that passes through the point $s$ at the initial time $t=0$ is given by $x = s + u_{0}(s)t$, and the velocity is constant, $u = u_{0}(s)$,  along the characteristic. In addition, it easily follows that along the characteristic
\[ \frac{\partial u}{\partial x}(x,t) = \frac{u_{0}'(s)}{1 + u_{0}'(s)t}\] and
\[ \rho(x,t) = \frac{\rho_{0}(s)}{1 + u_{0}'(s)t}.\]
This implies that 
\[\lim_{t \rightarrow -u_{0}'(s)^{-1}}  (\rho,u) = (+\infty,+\infty)\]
along the characteristic. Hence, we can only define a smooth solution up until we approach the finite-blowup time $t = -u_{0}'(s)^{-1}$ where the gradient and density exhibit singularities.

\subsection{Distribution solutions}
Define a two-dimensional weighted delta function $w(s)\delta_{\mathcal{C}}$ supported on a smooth curve $\mathcal{C}$ parametrized by $t=t(s)$, $x=x(s)$ $(a<s<b)$ by
\[ <w(s)\delta_{\mathcal{C}}, \psi(x,t)> = \int_{a}^{b} w(s)\psi(x(s),t(s))\sqrt{x'^{2} + t'^{2}}\,ds \]
for all $\psi \in C_{0}^{\infty}(\mathbb R_{+}^{2})$. 

Then we define two distributions in the dual space $(C_{0}^{\infty}(\mathbb R_{+}^{2}))^{*}$ by
\begin{equation}\label{distribution1}
<h, \phi> = \int_{0}^{\infty} \int_{-\infty}^{\infty} h\phi \,dxdt + <w(s)\delta_{\mathcal{C}},\phi>,
\end{equation}
and
\begin{equation}\label{distribution2}
<h u, \phi> = \int_{0}^{\infty} \int_{-\infty}^{\infty} hu\phi \,dxdt + <w(s)\delta_{\mathcal{C}},\phi>, 
\end{equation}
for $\phi \in C_{0}^{\infty}(\mathbb R_{+}^{2})$ and $h$ locally integrable.

\begin{definition}
Define two distributions $\in(C_{0}^{\infty}(\mathbb R_{+}^{2}))^{*}$ as in (\ref{distribution1}) and (\ref{distribution2}). Then we say that $\rho$ is a (Dirac delta) distribution solution to (\ref{onedim1}) if
\[ <\rho,\phi_{t}> + <\rho u,\phi_{x}> = 0 \]
for all $\phi \in C_{0}^{\infty}(\mathbb R_{+}^{2})$. For the Burgers equation, its standard definition of a generalized or weak solution will be used.
\end{definition}

\subsection{The Riemann problem}
Consider (\ref{onedim1})-(\ref{onedim2}) with the following initial condition

\begin{equation} 
(\rho_{0}(x),u_{0}(x)) =
\left\{ 
\begin{array}{c l}  
  (\rho_{l},u_{l}), & x < 0\\
  (\rho_{r},u_{r}), & x > 0
\end{array}\right.
\end{equation}

where $\rho_{l}$, $\rho_{r} \geq 0$.
This is nonstrictly hyperbolic with a repeated eigenvalue $\lambda = u$ and corresponding eigenvector $\vec{r} = (1,0)^{T}$. The eigenvalues are linearly degenerate, thus the solutions can be constructed by contact discontinuities, a vacuum or delta shock wave connecting the two constant states $(\rho_{l},u_{l})$, $(\rho_{r},u_{r})$.

\textbf{case(i)} \ $u_{l} < u_{r}$ 

A vacuum occurs between two contact discontinuities. A self-similarity transform $\psi = x/t$ yields the solution
\begin{equation} 
(\rho,u)(x,t) =
\left\{
\begin{array}{c l} 
  (\rho_{l},u_{l}), & x/t < u_{l}\\
  ( 0 , x/t ), & u_{l} \leq x/t \leq u_{r}\\
  (\rho_{r},u_{r}), & x/t > u_{r}.
\end{array}
\right.
\end{equation}

\textbf{case(ii)} \ $u^{*} := u_{l} = u_{r}$ 

It easily follows that the left and right constant states are connected by a contact discontinuity along $x = u^{*}t$.

\textbf{case(iii)} \ $u_{l} > u_{r}$ 

This is the case in which a delta-shock formation occurs i.e. the density $\rho$ is the sum of an $L^{\infty}$ function and a weighted delta function along the shock curve.
\begin{equation} \label{deltashock}
(\rho,u)(x,t) =
\left\{
\begin{array}{c l}
  (\rho_{l},u_{l}),                          & x/t < \sigma\\
  ( w(t)\delta(x - \sigma t) , u_{\delta} ), & x/t = \sigma\\
  (\rho_{r},u_{r}),                          & x/t > \sigma,
\end{array}
\right.
\end{equation}

where \[\sigma = u_{\delta} = \frac{1}{2}(u_{l} + u_{r}) \text{,} \; u_{l} > \sigma > u_{r} \text{,}\; \text{ and } \; w(t) = \frac{t}{\sqrt{1 + \sigma^{2}}}([\rho u]-\sigma[\rho]).\]
Note that $u$ is given the value $u_{\delta}$ along the shock. This is justified because if the values of $u$ are changed on a set of 2-d Lebesgue measure zero, then $u$ would still remain a weak solution for Burgers equation. 

\begin{figure}
\begin{minipage}[b]{0.45\linewidth} 
\centering
\includegraphics[width=5cm]{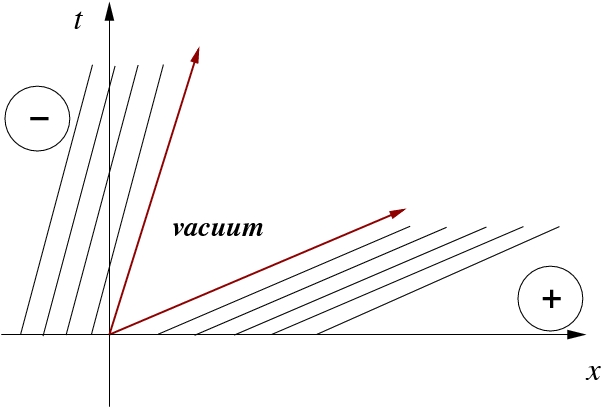}
\caption{Case $u_{l} < u_{r}$}
\end{minipage}
\hspace{0.5cm} 
\begin{minipage}[b]{0.45\linewidth}
\centering
\includegraphics[width=5cm]{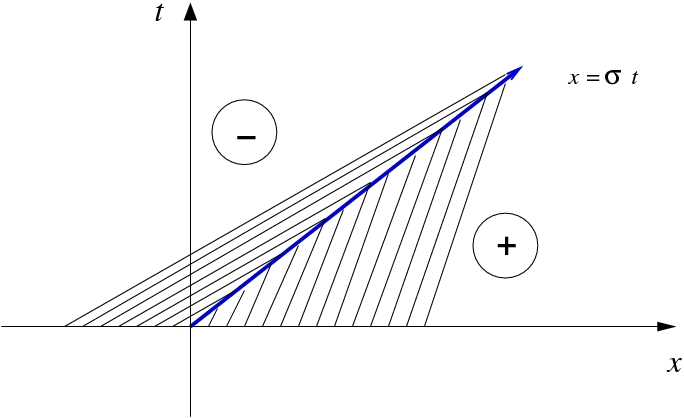}
\caption{Case $u_{l} > u_{r}$}
\end{minipage}
\end{figure}


\section{Filters and averaging kernels}\label{Filters}
In this section, we define the acceptable filters used in averaging the chosen quantities in the equations. For a given real-valued function $f$, we introduce an averaging kernel $g$ and define the filtering of $f$ by the convolution operator
\[ \overline{f} = g*f = \int_{\mathbb R^{n}} g(x-y)f(y)\,dy \]

for $x$ belonging to some given domain $\Omega \subseteq \mathbb R^{n}$.\\
We shall only consider averaging kernels $g \in W^{1,1}(\mathbb R^{n})$ that satisfy the following:

\begin{enumerate}
\item \[ \int_{\mathbb R^{n}} g \,d\ob{x} = 1 \] 
\item \[ g(\ob{x}) > 0 \text{ for all } \ob{x} \in \Omega \]
\item \[ |\ob{x}_{1}| \leq |\ob{x}_{2}| \Rightarrow |g(\ob{x}_{1})| \leq |g(\ob{x}_{2})| \]
\item \[ |\ob{x}_{1}| = |\ob{x}_{2}| \Rightarrow |g(\ob{x}_{1})| = |g(\ob{x}_{2})| \]
\item \[ \lim_{|\ob{k}|\rightarrow \infty}\ob{k} \hat{g}(\ob{k}) = 0,  \text{ where } \hat{g} \text{ is the Fourier transform of } g \]
\end{enumerate}
In the physical sense, the averaging should only provide non-negative weight to particles, have no preferential direction, and should give more weight to particles that are physically closer. In addition, we prescribe a parameter, $\alpha > 0$, to such a filter such that
\[ g^{\alpha} = \frac{1}{\alpha} g\left(\frac{\ob{x}}{\alpha}\right). \]
This parameter $\alpha$ acts as a scaling of the kernel and controls the level of filtering. One example of a commonly studied filter is the Helmholtz filter $f = \overline{f} - \alpha^{2} \Delta \overline{f}$ corresponding to the averaging kernel
\[g(x) = \frac{1}{2\alpha} exp\left(-\frac{|x|}{\alpha}\right).\]

\section{Flux Averaging and the Method of Observable Divergence}\label{Averaging}
This section will examine how the particular spatial averaging method is applied to the transport equation including background on its development. 

For the Burgers equation, the filtering of the convective velocity i.e. the convectively filtered Burgers (CFB) equation $u_{t} + \overline{u}u_{x} = 0$, provided a valid regularization by preventing the crossing of characteristics \cite{BhatHS:06a, Mohseni:08w}. With the success of this regularization, an extension to this averaging was sought for more general conservation laws. One such technique--dubbed the observable divergence for a vector field--entails an averaging of the flux term \cite{Mohseni:09w} . More precisely, if we are given a scalar function $f$ and a vector field $\textbf{V}$, we define the observable divergence of the vector field $\ob{F} = f\ob{V}$ by
\begin{equation}\label{odiv}
odiv (F) = \overline{f} \ \nabla \cdot \ob{V} + \overline{\ob{V}} \cdot \nabla f.
\end{equation}
Then conservation laws of the type $f_{t} + \nabla \cdot (f \ob{V}) = 0$ are modified into its `observable' conservation law $f_{t} + odiv(F) = 0$. Furthermore, Mohseni later derived the `observable' conservation law equations from basic principles \cite{Mohseni:09w}. 

To explain how the CFB equation inspired this particular averaging, note that the flux term for the one-dimensional Burgers equation, $(\frac{u^{2}}{2})_{x} = (\frac{uu}{2})_{x}$, may be expanded as 
\begin{equation}\label{split}
\frac{1}{2}(uu_{x} + u_{x}u). 
\end{equation}
With this in mind, Norgard and Mohseni \cite{Mohseni:09j} developed the observable divergence method as a generalization of the averaging used in the CFB equation by observing that averaging the non-differentiated terms in (\ref{split}) 
\begin{equation}\label{splitav}
\frac{1}{2}(\overline{u}u_{x} + u_{x}\overline{u})
\end{equation}
yields the relationship $odiv(\frac{uu}{2}) = \frac{1}{2}(\overline{u}u_{x} + u_{x}\overline{u}) = \overline{u}u_{x}$. Therefore, we may formally obtain the CFB equation by applying the method of observable divergence to Burgers equation i.e.
\[ u_{t} + odiv\left(\frac{u^{2}}{2}\right) = 0 \ \Longleftrightarrow \ u_{t} + \overline{u}u_{x} = 0.  \] 
We extend this averaging for Burgers equation in multiple-dimensions by
\[ \ob{u}_{t} + \overline{\ob{u}} \cdot \nabla \ob{u} = 0. \] 
Clearly, (\ref{odiv}) reasonably follows from this concept of expanding the flux and filtering the non-differentiated terms as illustrated in (\ref{split}) and (\ref{splitav}). Now consider the continuity equation $\rho_{t} + \nabla \cdot{(\rho \ob{u})} = 0$ and apply the observable divergence method to
$\ob{F} = (\rho \ob{u})$ to obtain the equation 
\[\rho_{t} + \overline{\rho} \ \nabla \cdot \ob{u} + \overline{\ob{u}} \cdot \nabla \rho = 0.\]
We shall call this modified continuity equation the observable continuity equation. Hence, the observable transport equations--the system combining the observable continuity equation and CFB equation--is derived by applying the method of observable divergence to the transport equations.

\section{Observable equations}\label{Observable equations}
Consider the observable transport equations
\begin{equation}\label{Equation2}
  \left\{\begin{array}{r@{\;=\;}l@{\quad}l}
      	\rho_{t} + \overline{\rho} \nabla \cdot \ob{u} + \overline{\ob{u}} \cdot \nabla \rho & 0 \\
 	 			\ob{u}_{t} + \overline{\ob{u}} \cdot \nabla \ob{u} & 0. \\
        \end{array}
 \right.
 \end{equation}
 
This section will address the very importannt and desirable property concerning the existence of solutions to the observable transport equations in some given domain $\Omega \times I \subseteq \mathbb R^{n} \times \mathbb R$ and satisfying some given initial data
\begin{center}
$(\rho(\ob{x},0),\ob{u}(\ob{x},0)) = (\rho_{0}(\ob{x}),\ob{u}_{0}(\ob{x}))$.
\end{center}
The results of this section is summarized as follows. An existence and uniqueness result for the multi-dimensional transport equations with smooth initial conditions will be established. First, the case for the CFB equation is covered. Next, a fixed point theorem is used to show the existence and uniqueness of a global broad solution for the observable continuity equation. Then it is shown that smooth initial conditions will yield classical solutions. The use of this spatial averaging as a legitimate alternative to vanishing viscosity requires it to correctly model the delta-shock phenomena. Therefore the Riemann problem for the one-dimensional problem is considered since it provides the simplest example of discontinuous initial data that provide solutions that contains a delta-shock distribution solution. It will be shown that the observable transport equations share the exact delta-shock distribution solution as the non-regularized transport equations.
 
\subsection{Existence and Uniqueness} 
If we start with smooth initial conditions 
$u_{0}(x) \in C^{1}(\mathbb R^{n})$, we can show the existence and uniqueness of the velocity field $u(x,t) \in C^{1}(\mathbb R^{n}, \mathbb R)$ to the initial value problem
\begin{equation}\label{Equation3}
  \left\{\begin{array}{r@{\;=\;}l@{\quad}l}
      	\ob{u}_{t} + \overline{\ob{u}} \cdot \nabla \ob{u} & 0 \\
 	 			\ob{u}(x,0) & \ob{u}_{0}(\ob{x}). \\
        \end{array}
 \right.
 \end{equation}
We briefly describe the procedure to showing this, but the reader is referred to \cite{Mohseni:08w} for the complete proof. Associate the map $\phi_{t}(\ob{s}): \mathbb R^{n} \mapsto \mathbb R^{n}$ satisfying
\[\frac{\partial }{\partial t} \phi_{t}(\ob{s}) = \overline{\ob{u}} (\phi_{t}(\ob{s}),t)\]
with initial condition $\phi_{t=0}(\ob{s}) = \ob{s}$. This map can be interpreted as mapping the initial position to its position at time $t$ and is called the Lagrangian map. Thus (\ref{Equation3}) yields that $\dfrac{d}{dt} \ob{u}(\phi_{t}(\ob{s}),t) = 0$, which implies that
\[u(\phi_{t}(\ob{s}),t) = u(\phi_{t=0}(\ob{s}),0) = u_{0}(\ob{s}).\]
One can show that the Jacobian of $\phi_{t}$, $J(\phi_{t})$, is nonzero for all $\ob{s}\in \mathbb R^{n}$ and for every $t\geq 0$. Hence, the mapping $t \mapsto \phi_{t}(\ob{s})$ is a diffeomorphism for each fixed $t$. From this result, it follows that the velocity field is given by 
\[ \ob{u}(\ob{x},t) = \ob{u}_{0}(\phi^{-1}_{t}(\ob{x})). \]
\textit{Remark}: Similarly, if the initial condition $\ob{u}_{0}(\ob{x}) \in L^{\infty}(\mathbb R^{n})$, then there exists a unique solution $\ob{u}(\ob{x},t) \in L^{\infty}(\mathbb R^{n} ,\mathbb R)$.

Now we are only left with determining a solution for the density function, $\rho(x,t)$ to the initial value problem
\begin{equation}\label{Equation4}
  \left\{\begin{array}{r@{\;=\;}l@{\quad}l}
      	\rho_{t} + \overline{\rho} \ \nabla \cdot \ob{u} + \overline{\ob{u}} \cdot \nabla \rho & 0 \\
 	 			\rho(x,0) & \rho_{0}(x). \\
        \end{array}
 \right.
 \end{equation}
Initially, we prove the global existence of a bounded measurable solution that satisfies the observable continuity equation along the characteristic curves. From there, we show that if we choose the initial condition to be continuously differentiable and impose further regularity conditions on the velocity, then the previous bounded measurable solution turns out to be continuously differentiable as well.

\begin{definition}
A closed region $D = \Omega \times I \subseteq \mathbb R^{n} \times \mathbb R$ is called a domain of determinacy for (\ref{Equation4}) provided that the following holds. For every $(\xi, \tau) \in D$, the characteristic curve $\{(x(t;\xi,\tau),t)\}$ is entirely contained inside $D$. 
\end{definition}

From above, recall the characteristics are defined to be the map $\phi_{t}(s)$, and we showed that these characteristic curves are nonintersecting for different initial point $s$. Moreover, it was shown that the velocity was constant along these curves with respect to time. From this point on, we shift our perspective slightly. From the implicit function theorem, we can find a $C^{1}$ mapping $g:\mathbb R \rightarrow \mathbb R^{n}$ such that $\phi_{t}(s) = \phi_{t}(g(t)) =: x(t)$ for $t\geq 0$. That is the characteristics may be parametrized on our domain of choice by the curve $(x(t),t)$. As a consequence, for every point $(\xi,\tau)$ in our domain, we denote $x(t;\xi,\tau)$ to be the unique characteristic curve passing through $(\xi,\tau)$ (see Figure 3 below). Before stating and proving our desired result, let us first introduce several concepts and definitions.

\begin{figure}
\begin{center}
\includegraphics[width=0.9\linewidth]{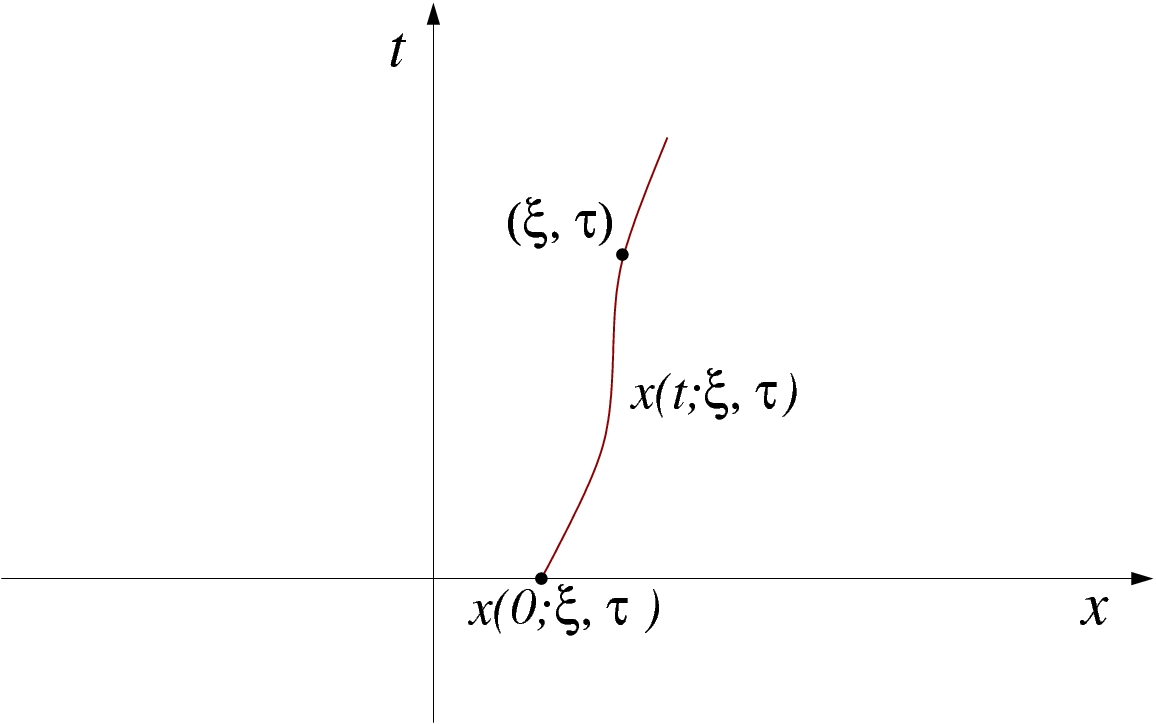}
\caption{\scriptsize Characteristic passing through $(\xi,\tau)$ }
\end{center}
\vskip -0.00 true in
\end{figure}

So suppose that $D$ is a domain of determinacy, and let $X$ be the Banach space of bounded measurable functions $v:D \mapsto \mathbb R$ with respect to the weighted essential supremum-norm
\[ \|v\| = \textrm{ess} \sup_{D} e^{-2Lt} |v(x,t)| \]
where $L$ is a fixed bound for $\nabla \cdot u$ on $D$.  

Furthermore, the characteristic curves $x(t;\xi,\tau)$ clearly satisfies $\dfrac{dx}{dt} = \overline{u}(x(t),t)$. Thus
\[ \dfrac{d\rho}{dt}(x(t;\xi,\tau),t) = -\overline{\rho}(x(t;\xi,\tau),t) \ \nabla\cdot u(x(t;\xi,\tau),t). \]

This observation motivates our definition of a broad solution to (\ref{Equation4}).

\begin{definition}
Let $H:X \rightarrow X$ be the map $H[\rho](x,t) = -\overline{\rho}(x,t) \nabla \cdot u(x,t)$. We say an integrable function $\rho:D \rightarrow \mathbb R$ is a \textit{broad} 
solution for the problem (\ref{Equation4}) provided that 
\[ \rho(\xi,\tau) = \rho_{0}(x(0;\xi,\tau)) + \int_{0}^{\tau} H(x(t;\xi,\tau),t,\rho(x(t;\xi,\tau),t)) \,dt \]
holds for $(\xi,\tau) \in D$ almost everywhere.
\end{definition}
Now the following lemma will prove useful in establishing a Lipschitz-type inequality required in Theorem 1 where a fixed point theorem is used to show the existence and uniqueness of a global broad solution.

\begin{lemma}
Consider again the map $H[\rho](x,t)$ as defined previously. Then
\[ \left|H[\rho_{1}](x,t) - H[\rho_{2}](x,t)\right| \leq Le^{2Lt} \|\rho_{1} - \rho_{2}\| \]
for $(x,t)\in D$ and $\rho_{1}$,$\rho_{2} \in X$. 
\end{lemma}
\begin{proof}
\begin{align*}
e^{-2Lt}|H[\rho_{1}](x,t) - H[\rho_{2}](x,t)| \leq {} & L e^{-2Lt}| \overline{\rho}_{1}(x,t) - \overline{\rho}_{2}(x,t)| \\
\leq {} & L \ \textrm{ess}\sup_{D} e^{-2Lt}|\overline{\rho}_{1}(x,t) - \overline{\rho}_{2}(x,t)| \\
\leq {} & L\|\overline{\rho}_{1}-\overline{\rho}_{2}\| \leq L\| \rho_{1} - \rho_{2}\|. \qedhere
\end{align*}
\end{proof}

\begin{theorem}
Let $D = \Omega \times I$ be our domain of determinacy. Suppose that $\rho_{0}:\Omega \rightarrow \mathbb R$ is a bounded measurable function and $u_{0} \in W^{1,\infty}(\Omega)$. Then there exists a unique bounded measurable solution $u:D \rightarrow \mathbb R$ to the initial value problem (\ref{Equation3}) and a unique broad solution $\rho:D \rightarrow \mathbb R$ to the initial value problem (\ref{Equation4}).
\end{theorem}
\begin{proof}

Define the integral operator $T:X \rightarrow X$ to be
$$T\rho(\xi,\tau) = \rho_{0}(x(0;\xi,\tau)) + \int_{0}^{\tau} H(x(t;\xi,\tau),t,\rho(x(t;\xi,\tau),t)) \,dt$$
\textbf{Claim:} The map $T:X \rightarrow X$ is a contraction map.

For $(x,t)\in D$, $\rho_{1}$,$\rho_{2} \in X$, and by Lemma 1,
\begin{align*}
|T\rho_{1}(\xi,\tau) - T\rho_{2}(\xi,\tau)|
\leq {} & \int_{0}^{\tau}|H[\rho_{1}(x(t;\xi,\tau),t)]-H[\rho_{2}(x(t;\xi,\tau),t)]|\,dt \\
\leq {} & \int_{0}^{\tau} Le^{2Lt} \|\rho_{1} - \rho_{2}\| \,dt \leq L \|\rho_{1} - \rho_{2}\| \int_{0}^{\tau} e^{2Lt} \,dt \\
\leq {} & \|\rho_{1} - \rho_{2}\| \left(\frac{1}{2} e^{2Lt}\Big |_{0}^{\tau}\right) \leq \frac{1}{2}\|\rho_{1} - \rho_{2}\|(e^{2L\tau} - 1) \\
\leq {} & \frac{1}{2} e^{2L\tau} \|\rho_{1} - \rho_{2}\|
\end{align*}

It follows that 
\[ e^{-2L\tau} |T\rho_{1}-T\rho_{2}| < \frac{1}{2}\|\rho_{1} - \rho_{2}\|. \] 
That is, 
\[ \|T\rho_{1} - T\rho_{2}\| < \frac{1}{2} \|\rho_{1} - \rho_{2}\|. \]

Hence, by the Banach fixed point theorem \cite{Evans:98a}, there exists a unique element $\rho \in X$ such that $T\rho = \rho$. The existence and uniqueness of $u$ has been discussed previously.\qedhere
\end{proof}

\subsubsection*{Classical Solutions}
Let us now prove that if we have $C^{1}$ initial data, then the unique bounded measurable solutions to (\ref{Equation2}) are in fact within the class of $C^{1}$ functions.

\begin{theorem}
Let $\rho_{0}:\Omega \rightarrow \mathbb R$ and $u_{0}:\Omega \rightarrow \mathbb R$ be continuously differentiable and assume $D=\Omega \times I$ is a compact domain of determinacy. Moreover, suppose that the second partial derivatives of $\textbf{u}_{0}$ are bounded on $\Omega$. Then there exist unique solutions $u(x,t)$, $\rho(x,t) \in C^{1}(D)$ to the IVP (\ref{Equation2}).
\end{theorem}

\begin{proof}
From the previous theorem, $\rho$ is the uniform limit of the sequence $\rho^{(m)}$, for $m\geq 1$ defined inductively by
\begin{center}
$\rho^{(m+1)} = T\rho^{(m)}$ with $\rho^{(0)}(x,t) = \rho_{0}(x)$,
\end{center}
or equivalently as the sequence
\[ \rho^{(m+1)}(\xi,\tau) = \rho_{0}(x(0;\xi,\tau)) - \int_{0}^{\tau} (\nabla \cdot u) \overline{\rho}^{(m)} \,dt. \]

Since the initial data $\rho_{0}$ is continuously differentiable, we deduce inductively that if $\rho^{(m)}$ is $C^{1}$ then $\rho^{(m+1)}$ is $C^{1}$ for each $n \geq 0$. 
Since the uniform limit of a sequence of continuously differentiable functions with respect to the supremum norm is at least continuous, we are only left to show that $\rho$ has continuous partial derivatives on $D$.

Take the partial derivative of this with respect to $\xi_{i}$ for $1 \leq i \leq n$ to get
\[ \frac{\partial}{\partial \xi_{i}} \rho^{(m+1)}(\xi,\tau) = \sum_{k=1}^{n}\frac{\partial \rho_{0}}{\partial x_{k}} \frac{\partial x_{k}(0;\xi,\tau)}{\partial \xi_{i}} - \int_{0}^{\tau} \overline{\rho}^{(m)} 
\frac{\partial}{\partial \xi_{i}} (\nabla \cdot u) + \sum_{k=1}^{n} \nabla \cdot u \frac{\partial \overline{\rho}^{(m)}}{\partial x_{k}} \frac{\partial x_{k}(t; \xi,\tau)}{\partial \xi_{i}} \,dt \]

where the integrals are evaluated along the curve $t \mapsto x(t;\xi,\tau)$. Now for each integer $m\geq 1$ and every $\tau \geq 0$, define
\[ Y_{m}(\tau) = \sup_{(\xi,\tau)} \left| \frac{\partial \rho^{(m)}}{\partial \xi_{i}}(\xi,\tau) \right| \] and
\[ Z_{m}(\tau) = \sup_{(\xi,\tau)} \left| \frac{\partial \rho^{(m+1)}}{\partial \xi_{i}}(\xi,\tau) - \frac{\partial \rho^{(m)}}{\partial \xi_{i}}(\xi,\tau) \right| \]

Then there are constants $C_{1}$, $C_{2}$ such that
\[ Y_{m+1}(\tau) \leq C_{1} + C_{2}\int_{0} ^{\tau} Y_{m}(t)\, dt. \]

This inequality implies that $Y_{m}(\tau) \leq C_{1}e^{C_{2}\tau}$. Therefore $\{Y_{m}\}$ is uniformly bounded.

Since $T$ was shown to be a contraction map, we can find a constant $C_{3}$ such that 
\begin{center}
$|\rho^{m+1}(\xi,\tau) - \rho^{m}(\xi,\tau)| \leq C_{3}2^{-m}$ for all $(\xi,\tau) \in D$, $m \geq 0$.
\end{center}

With this we have, for some constant $C_{4}$, the following inequality 
\begin{center}
$Z_{m+1} \leq 2^{-m}C_{3} + C_{4}\int_{0} ^{\tau} Z_{m}(t)\,dt$.
\end{center}
This inequality implies that $Z_{m} \leq 2^{-m}(C_{3}+1)e^{2C_{4}\tau}$. Hence, $\sum_{m=0}^{\infty} Z_{m}$ converges uniformly. It follows the partial derivatives 
$\frac{\partial \rho^{(m)}}{\partial x_{i}} \rightarrow \frac{\partial \rho}{\partial x_{i}}$ uniformly therefore $\frac{\partial \rho}{\partial x_{i}}$ is continuous on $D$. 
A similar argument can be applied to the case concerning $\frac{\partial \rho}{\partial t}$. As briefly discussed above, the result for $u$ has been covered above.\qedhere
\end{proof}

\subsection{The Riemann Problem}
The existence of solutions to the observable transport equations has been covered and it was revealed that the smoothness of the solutions depend on the smoothness of the initial data, but it has yet to be verified if the regularized equations recover the delta-shock solutions, as is normally accomplished with vanishing viscosity. The basic Riemann problem is a well-known example of a discontinuous, piecewise constant initial data that provides a Dirac delta-shock solution. Remarkably, it will be shown that for arbitrary $\alpha > 0$, the observable equations may contain delta-shock distribution solutions as well. In fact, it shares the same distribution solution as the non-regularized equations. So consider the Riemann problem for our regularized system
\begin{equation}\label{onedimreg}
  \begin{array}{r@{\;=\;}l@{\quad}l}
      	\rho_{t} + \overline{u}\rho_{x} + \overline{\rho} u_{x} & 0 \\
 	 			u_{t} + \overline{u}u_{x} & 0 \\
 	 			(\rho(x,0),u(x,0)) & (\rho_{0},u_{0})
        \end{array}
\end{equation}

where 
\begin{equation}
(\rho_{0}(x),u_{0}(x)) =
\left\{ 
\begin{array}{c l}  
  (\rho_{l},u_{l}), & x < 0\\
  (\rho_{r},u_{r}), & x > 0.
\end{array}\right.
\end{equation}

We shall restrict ourselves to the case $u_{l} > u_{r}$ since this is the case where the solution is a delta-shock wave for (\ref{onedim1})-(\ref{onedim2}). 
Recall the characteristic curves via the Lagrangian particle map $\phi_{t}(s)$. Again, we have
\begin{equation}\label{char}
\frac{d\phi_{t}(s)}{dt} = \overline{u}(\phi_{t}(s),t).
\end{equation}

The solutions for the Riemann problem for CFB equation is studied in \cite{BhatHS:09a}. Here, the CFB Riemann problem is solved using the method of characteristics instead of using the usual idea of weak solutions involving test functions and integration by parts. We cover the necessary results from \cite{BhatHS:09a}. It is shown that
\begin{equation} 
u(x,t) =
\left\{
\begin{array}{c l}
  u_{l}, & x/t < \sigma\\
  u_{r}, & x/t > \sigma
\end{array}
\right.
\end{equation}
where $\sigma = \dfrac{u_{l} + u_{r}}{2}$.

\textit{Remark:} This solution is independent of $\alpha$.

Furthermore, without any restriction on $u_{l}$ and $u_{r}$, it follows that (\ref{char}) implies
\begin{equation} \label{gensoln}
\overline{u}(x,t) =
\left\{
\begin{array}{c l}
  u_{l}+(u_{r}-u_{l})\varphi^{-}(\frac{x-\sigma}{\alpha}),  & x/t < \sigma\\
  u_{\delta},                                               & x/t = \sigma\\
  u_{r}+(u_{r}-u_{l})\varphi^{+}(\frac{x-\sigma}{\alpha}),  & x/t > \sigma.
\end{array}
\right.
\end{equation}
Here $\varphi^{-}(x) = \int_{-\infty}^{x} g^{-}(y)\,dy$ and $\varphi^{+}(x) = \int_{x}^{\infty} g^{+}(y)\,dy$ where the averaging kernel is expressed as a piecewise function
\begin{equation}
g(x) =  
\left\{
\begin{array}{c l}
  g^{-}(x), & x < 0 \\
  g^{+}(x), & x > 0.
\end{array}
\right.
\end{equation}

Due to conditions placed on the averaging kernel, $\varphi^{\pm}$ decay to zero as $|x| \longrightarrow \infty$. Therefore it follows from (\ref{gensoln}) that
\begin{equation} 
\lim_{\alpha \rightarrow 0^{+}} \overline{u}(x,t) = 
\left\{
\begin{array}{c l}
  u_{l},   & x/t < \sigma \\
  u_{\delta},  & x/t = \sigma \\
  u_{r},   & x/t > \sigma.
\end{array}
\right.
\end{equation}
In the case $u_{l} > u_{r}$, this limit recovers the same solution for $u(x,t)$ as that in (\ref{deltashock}).

By the method of characteristics, the correct solution for $u$ is obtained, however, we are left to show that the delta-shock solution is preserved. So unlike the case for the CFB equation, we must express $\rho_{t} + \overline{u}\rho_{x} + \overline{\rho} u_{x} = 0$ in conservation form in order to define a distribution solution as in Section \ref{Preliminaries}. For brevity, let us restrict our attention to the Helmholtz filter, $u = \overline{u} - \alpha^{2}\overline{u}_{xx}$. Choose a test function $\phi$ with compact support and observe that
\begin{align*}
\rho_{t} + \overline{\rho}u_{x} + \overline{u}\rho_{x} = {} & \rho_{t} + (\overline{\rho}u + \rho \overline{u})_{x} - \overline{\rho}_{x}u - \rho\overline{u}_{x} \\
= {} & \rho_{t} + (\overline{\rho}u + \rho\overline{u})_{x} - \overline{\rho}_{x}(\overline{u} - \alpha^{2}\overline{u}_{xx}) - (\overline{\rho} - \alpha^{2}\overline{\rho}_{xx})\overline{u}_{x} \\
= {} & \rho_{t} + (\overline{\rho}u + \rho\overline{u})_{x} - (\overline{\rho}\hspace{0.5mm}\overline{u})_{x} + \alpha^{2}(\overline{\rho}_{x}\overline{u}_{x})_{x}.
\end{align*}
Integrating this against $\phi$ and using integration by parts gives us
\[ \iint_{R^{2}_{+}} \rho\phi_{x} + \rho\overline{u}\phi_{x} + \overline{\rho}(u-\overline{u})\phi_{x} + \alpha^{2}(\overline{\rho}_{x}\overline{u}_{x})\phi_x \,dxdt. \]
In light of this observation and in accordance with Section \ref{Preliminaries}, we say that $\rho$ is a distribution solution to the observable continuity equation if 
\[ <\rho,\phi_{t}> + <\rho \overline{u},\phi_{x}> + <\overline{\rho}(u - \overline{u}),\phi_{x}> + <\alpha^{2}\overline{\rho}_{x}\overline{u}_{x},\phi_{x}> = {} 0. \] 
The next theorem will show that the original delta-shock distribution solution to the transport equation also satisfies the regularized equations (\ref{onedimreg}) as well.

\begin{theorem} 
Let $\alpha > 0$ and suppose that $u_{l} > u_{r}$, $\sigma = \frac{u_{l}+u_{r}}{2} \in (u_{r},u_{l})$, and $[\rho u]-\sigma[\rho] > 0$. Furthermore, let $g$ be the Helmholtz filter. Then $(\rho,u)$ given by
\begin{equation}
\rho = h + w(s)\delta_{\mathcal{C}},
\end{equation}
where
\begin{equation} 
h(x,t) =
\left\{
\begin{array}{c l}
  \rho_{l}, & x/t < \sigma\\ 
  \rho_{r}, & x/t > \sigma,
\end{array}
\right.
w(s) = \frac{s}{\sqrt{1 + \sigma^{2}}}([\rho u]-\sigma[\rho]),
\end{equation}
\[ \mathcal{C}:= \{t = s, x = \sigma s: \ 0\leq s < \infty\}, \]
and
\begin{equation}
u(x,t) =
\left\{
\begin{array}{c l}
  u_{l},       & x/t < \sigma \\
  u_{\delta},  & x/t = \sigma \\
  u_{r},       & x/t > \sigma,
\end{array}
\right.
\end{equation}
satisfy the Riemann problem to (\ref{onedimreg}).
\end{theorem}

\begin{proof}
Recall that $u(x,t)$ and $\overline{u}(x,t)$ was determined using the method of characteristics. For the Helmholtz filter, it follows that
\begin{equation}
\overline{\rho}(x,t) = \overline{h} + w(s)\delta_{\mathcal{C}}
\left\{
\begin{array}{c l}
  \rho_{l} - \frac{1}{2}(\rho_{l}-\rho_{r})exp(\frac{(x-\sigma t)}{\alpha}),  & x/t < \sigma \\
  w(t)\delta(x-\sigma t),                                                     & x/t = \sigma \\
  \rho_{r} - \frac{1}{2}(\rho_{r}-\rho_{l})exp(-\frac{(x-\sigma t)}{\alpha}), & x/t > \sigma,
\end{array}
\right.
\end{equation}

\begin{equation} 
\overline{u}(x,t) = 
\left\{
\begin{array}{c l}
  u_{l} - \frac{1}{2}(u_{l}-u_{r})exp(\frac{(x-\sigma t)}{\alpha}),  & x/t < \sigma \\
  u_{\delta},                                                        & x/t = \sigma \\
  u_{r} - \frac{1}{2}(u_{r}-u_{l})exp(-\frac{(x-\sigma t)}{\alpha}), & x/t > \sigma,
\end{array}
\right.
\end{equation}

\begin{equation}
\overline{\rho}_{x}(x,t) =
\left\{
\begin{array}{c l}
  -\frac{1}{2\alpha}(\rho_{l}-\rho_{r})exp(\frac{(x-\sigma t)}{\alpha}),  & x/t < \sigma \\
  \frac{1}{2\alpha}(\rho_{r}-\rho_{l})exp(-\frac{(x-\sigma t)}{\alpha}),  & x/t > \sigma,
\end{array}
\right.
\end{equation}
and
\begin{equation} 
\overline{u}_{x}(x,t) = 
\left\{
\begin{array}{c l}
  -\frac{1}{2\alpha}(u_{l}-u_{r})exp(\frac{(x-\sigma t)}{\alpha}),  & x/t < \sigma \\
  \frac{1}{2\alpha}(u_{r}-u_{l})exp(-\frac{(x-\sigma t)}{\alpha}),  & x/t > \sigma.
\end{array}
\right.
\end{equation}
To complete the proof, it is left to show the following claim holds true.
\[ \textbf{Claim:} <\rho,\phi_{t}> + <\rho \overline{u},\phi_{x}> + <\overline{\rho}(u - \overline{u}),\phi_{x}> + <\alpha^{2}\overline{\rho}_{x}\overline{u}_{x},\phi_{x}> = {} 0\]
for all $\phi \in C^{\infty}_{0}(\mathbb R^{2}_{+})$.

To prove this claim we compute each term on the right-hand side then verify that their sum is zero.
\begin{align*}
\text{(i)} <\rho,\phi_{t}> + <\rho \overline{u},\phi_{x}>
= {} & \iint_{\mathbb R_{+}^{2}} \rho \phi_{t} + (\rho\bar{u})\phi_{x} \,dx dt \\
= {} & \int_{0}^{\infty}\int_{x/\sigma}^{\infty} \rho_{l}\phi_{t} \,dtdx + \int_{0}^{\infty}\int_{0}^{x/\sigma}\rho_{r}\phi_{t}\,dtdx \\
& + \int_{0}^{\infty} w(t)\left(\phi_{t}(\sigma t,t) + u_{\delta}\phi_{x}(\sigma t,t)\right)\sqrt{1 + \sigma^{2}} \,dt \\
& + \int_{0}^{\infty}\int_{-\infty}^{\sigma t} \rho_{l}(u_{l}-\frac{1}{2}(u_{l}-u_{r})exp\left(\frac{x-\sigma t}{\alpha}\right)\phi_{x}\,dxdt \\
& + \int_{0}^{\infty}\int_{\sigma t}^{\infty} \rho_{r}(u_{r}-\frac{1}{2}(u_{r}-u_{l})exp\left(\frac{-(x-\sigma t)}{\alpha}\right)\phi_{x}\,dxdt \\
= {} & -\int_{0}^{\infty}\int_{-\infty}^{\sigma t} \frac{1}{2}\rho_{l}(u_{l}-u_{r})exp\left(\frac{x-\sigma t}{\alpha}\right)\phi_{x}\,dxdt \\
& - \int_{0}^{\infty}\int_{-\infty}^{\sigma t} \frac{1}{2}\rho_{r}(u_{r}-u_{l})exp\left(\frac{-(x-\sigma t)}{\alpha}\right)\phi_{x}\,dxdt.
\end{align*}

Similar calculations for the other terms gives us
\begin{align*}
\text{(ii)} <\overline{\rho}(u-\overline{u}),\phi_{x}>
= {} & \int_{0}^{\infty}\int_{-\infty}^{\sigma t}\frac{1}{2}\rho_{l}(u_{l}-u_{r})exp\left(\frac{x-\sigma t}{\alpha}\right)\phi_{x} \,dxdt \\
& -\int_{0}^{\infty}\int_{-\infty}^{\sigma t} \frac{1}{4}(\rho_{l}-\rho_{r})(u_{l}-u_{r})exp\left(\frac{2(x-\sigma t)}{\alpha}\right)\phi_{x}\,dxdt\\
& +\int_{0}^{\infty}\int_{\sigma t}^{\infty} \frac{1}{2}\rho_{r}(u_{r}-u_{l})exp\left(-\frac{x-\sigma t}{\alpha}\right)\phi_{x}\,dxdt \\
& -\int_{0}^{\infty}\int_{\sigma t}^{\infty} \frac{1}{4}(\rho_{r}-\rho_{l})(u_{r}-u_{l})exp\left(-\frac{2(x-\sigma t)}{\alpha}\right)\phi_{x}\,dxdt
\end{align*}
and
\begin{align*}
\text{(iii)} <\alpha^{2}\overline{\rho}_{x}\overline{u}_{x},\phi_{x}>
= {} & \int_{0}^{\infty}\int_{-\infty}^{\sigma t} \frac{1}{4}(\rho_{l}-\rho_{r})(u_{l}-u_{r})exp\left(\frac{2(x-\sigma t)}{\alpha}\right)\phi_{x}\,dxdt\\
& \int_{0}^{\infty}\int_{\sigma t}^{\infty} \frac{1}{4}(\rho_{r}-\rho_{l})(u_{r}-u_{l})exp\left(-\frac{2(x-\sigma t)}{\alpha}\right)\phi_{x}\,dxdt.
\end{align*}
Hence, summing up these terms gives us zero i.e. (i) + (ii) + (iii) = 0\qedhere
\end{proof}

For completeness, let us briefly discuss the case for $u_{l} < u_{r}$. Recall that the original solution is a rarefaction wave with a vacuum and does not contain a delta function. The CFB equation will capture a non-entropic weak solution of Burgers rather than the rarefaction wave. This is remedied by replacing the Riemann initial condition with a filtered approximation. That is, introduce filtering to the initial condition independent from the CFB equation, then the CFB will recover the entropy solution for that filtered initial condition \cite{BhatHS:09a}. On the other hand, Norgard and Mohseni \cite{Mohseni:09h} showed the CFB equation recovers the entropy solution for a class of smooth initial conditions and proposed a more general program in recovering the entropy solution for discontinuous initial conditions. Such a program can prove useful in numerical simulations. More specifically, by spatially averaging Burgers equation via the CFB equation in conjunction with the discontinuous initial condition, the entropy solution is captured for the discontinuous initial data as filtering vanishes. As a caveat, this follows from \cite{Mohseni:09h} if one assumes that the CFB equation recovers the entropy solution for arbitrary continuous initial conditions. This behavior indicates how the initial smoothness can influence how effective spatial averaging will model our equations and suggests that the method of observable divergence may be more appropriate for more regular initial data.

\section{Concluding Remarks}
A spatial averaging technique--previously applied to the Burgers and Euler equations as a valid shock-regularization--was applied to the transport equations. Theory developed for the CFB equation combined with the past success of numerical analysis for the regularized equations was the primary motivator in studying the effects of averaging the transport equations. In particular, the method of observable divergence was introduced and applied to the transport equations. It was shown that utilizing this carefully chosen averaging technique will exhibit a unique classical solution for smooth initial conditions. But unlike the method of vanishing viscosity, the observable equations will preserve the original delta-shock solution for the Riemann problem regardless of the amount of filtering used. Nonetheless, the delta-shock is captured when the level of filtering vanishes. These results reinforce the validity that this technique is a legitimate alternative to the method of vanishing viscosity--further establishing its capability as an effective technique in properly modeling general conservation laws.

\section*{Acknowledgement} The research in this paper was partially supported by the AFOSR contract FA9550-05-1-0334. The authors would like to thank Greg Norgard for providing helpful discussions and his comments on a draft of this manuscript.

\bibliographystyle{unsrt}
\bibliography{transport_arxiv}

\end{document}